\documentclass[10pt]{article}

\usepackage{Variationen}
\usepackage{multicol}

\title{Veränderungen\\
\large über einen Satz von Timmesfeld\bigskip\\
I. Quadratic Actions}
\date{\today}
\author{Adrien Deloro}

\begin{document}

\maketitle

\liminaryquotation{Ich hätte glücklich geendet, aber diese Nro. 30, das Thema, riß mich unaufhaltsam fort. Die quartblätter dehnten sich plötzlich aus zu einem Riesenfolio, wo tausend Imitationen und Ausführungen jenes Themas geschrieben standen, die ich abspielen mußte.}

\abstract{We classify quadratic $\SL_2(\K)$- and $\sl_2(\K)$-modules by crude computation, generalizing in the first case a Theorem proved independently by F.G. Timmesfeld and S. Smith. The paper is the first of a series dealing with linearization results for abstract modules of algebraic groups and associated Lie rings.}

\renewcommand{\thefootnote}{}

\section*{General Foreword}


My hope is that the results gathered hereafter will suggest to the reader a simple idea: a whole chapter of representation theory could be written at the basic level of computations, 
without the help of algebraic geometry.
For I wish in this article and in others which may follow to study how much of geometric information on modules is already prescribed by the inner constraints of algebraic groups seen as abstract groups.
In a sense the problem is akin to the one solved by Borel and Tits in their celebrated work on abstract homomorphisms of algebraic groups. But here instead of morphisms between abstract groups we deal with abstract modules.

The central question is the following.
\begin{quote}
Let $\K$ be a field and $G$ be the abstract group of $\K$-points of an algebraic group. Is every $G$-module a $\K G$-module?
\end{quote}

It would be obscene to hope for a positive answer. The question is not asked literally; one should at least require the algebraic group to be reductive, if not simple. Moreover it may be necessary to bound the complexity of modules in some sense yet to be explained.

Since there is no $\K$-structure a priori and therefore no notion of a dimension over $\K$, one may focus on actions of finite nilpotence length, that is where unipotent subgroups act unipotently.
This setting seems to me more natural than that of $M_C$ modules (where centralizer chains are stationary); to support this impression one may bear in mind that the class of $M_C$ modules is not stable by going to a quotient, an operation which is likely to be relevant here. One could also make various model-theoretic assumptions, hoping that they would force configurations into the world of algebraic geometry; this did not seem natural either, since the first computations one can make are much too explicit for logic to play a deep role here. The future might bring contrary evidence; as for now, pure nilpotence seems more relevant.

So let us make our question more precise.

\begin{quote}
Let $\K$ be a field and $\bG$ be an algebraic group. Let $G = \bG_\K$. Understand the relationships between:
\begin{itemize}
\begin{multicols}{2}
\item
$\K G$-modules of finite length
\item
$G$-modules of finite length.
\end{multicols}
\end{itemize}
\end{quote}

One may be tempted to tackle the question by reducing it to actions of the associated Lie algebra. Two difficulties appear.
\begin{itemize}
\item
The Lie algebra ``remembers'' the base field, in a sense which we shall not explicit here; in any case one readily sees that actions of the Lie algebra can take place only in natural characteristic (that of the base field). Yet the group $G$ does not necessarily remember its base field, since there are various isomorphisms the extreme cases of which are over finite fields, such as $\SL_3(\F_2) \simeq \PSL_2(\F_7)$. These pathologies can be eliminated by reasonable assumptions on the algebraic group and on the field, and we shall deal only with decent cases.
\item
Above all the Lie algebra $\fg$ which is a $\K$-vector space, can appear only if $V$ is already equipped with a $\K$-linear structure; since it is not a priori, all one can hope for instead is the Lie ring, which is the (neither associative nor unitary) ring underlying the Lie algebra when one forgets its vector space structure. The representations of the Lie algebra are exactly the $\K \fg$-modules, where $\fg$ is seen as a Lie ring. Similarly, the universal object in this context will not be the enveloping algebra, but the enveloping ring.
\end{itemize}

And anyway nothing guarantees that reducing a group action to an action of its Lie ring (if possible) is any simpler than directly linearizing the module.
It thus looks like the introduction of the Lie algebra will not solve any question but bring new ones.
Our central problem extends as follows:

\begin{quote}
Let $\K$ be a field and $\bG$ be an algebraic group. Let $G = \bG_\K$ and $\fg = (\mathrm{Lie }\ \bG)_\K$ be its Lie algebra, seen as a Lie ring. Understand the relationships between:
\begin{itemize}
\begin{multicols}{2}
\item
$\K G$-modules of finite length
\item
$G$-modules of finite length
\item
$\K \fg$-modules of finite length
\item
$\fg$-modules of finite length.
\end{multicols}
\end{itemize}
\end{quote}

I want to speak for the idea that there are indeed good correspondences between these categories.
Here again this should not be taken literally: one may require the field to have sufficiently large characteristic and many roots of unity. 

The archetype of an effective linearization is the following result, proved independently by Stephen Smith and Franz Georg Timmesfeld (the latter mathematician actually did not require simplicity).
\begin{quote}
A simple $\SL_2(\K)$-module on which the unipotent subgroup acts quadratically is a $\K\SL_2(\K)$-module.
\end{quote}
It is not known whether the same holds over a skew-field. We shall not enter the topic, as all our results rely on heavy use of the Steinberg relations. Actually an alternative title might have been: ``$G$-modules and the Steinberg relations''.

The present work is constructed as a series of variations on the Smith-Timmesfeld theme, showing the unexpected robustness of the underlying computation. Encountered difficulties and provable results will provide equally important information: one must determine the limits of this computation in order to understand its deep meaning.
\begin{itemize}
\item
These variations will not be to the taste of geometers: from their point of view, I shall state only partial trivialities in an inadequate language.
But to defend these pages, and with a clear sense of proportions, 
I will appeal to the Borel-Tits famous result. The idea now is to understand to what extent inner constraints of abstract structures determine their representation theory. There is no rational structure here; everything is done elementarily.
\item
These variations could perhaps amuse group theorists, who must sometimes deal with austere objects with no categorical information.
Experts in finite group theory will nonetheless be upset by the lack of depth of my results, and by the efforts they cost: but the fields here may be infinite, and there is no character theory.
\item
The variations may at least be useful to logicians. Those with an interest in model-theoretic algebra often encounter abstract permutation groups; these sometimes turn out to be groups acting on abelian groups, and one needs results from more or less pure group theory to complete the discussion.
\end{itemize}

I confess that the present work takes place in a general context, far from model theory: I got carried away by the subject. 
To conclude this general foreword, I would love as much as the reader to suggest a conjecture describing in precise terms a phenomenon of ``linearity of abstract modules of structures of Lie type''; I would love to but I cannot, because it is too early.

My heart-felt thanks to Alexandre Borovik, for he believes in crazy ideas.

\section{The Setting}\label{S:setting}

In this article we study quadratic actions of $\SL_2(\K)$ and $\sl_2(\K)$ on an abelian group.\inmargin{Expanded.}

%
%

The articles cited in the next few paragraphs are by no means required in order to understand the hopefully self-contained present work. Only the reader with some knowledge of the topic will find interest in this introduction; the other reader may freely skip it. Such a liminary digression is merely meant to provide some historical background on the notion of quadraticity which lies at the center of our first article. The results we shall quote are not used anywhere and they bear no relationship to the rest of the series nor to its general spirit.

To the reader versed in finite group theory the word quadraticity will certainly evoke a line of thought initiated by J. Thompson: the classification of quadratic pairs, consisting of a finite group and a module with certain properties which we need not make precise. J. Thompson's seminal yet unpublished work \cite{Thompson} was quite systematically pursued by Ho \cite{Ho} among others, and more recently completed by A. Chermak \cite{Chermak} using the classification of the finite simple groups. This strain of results aims at pushing the group involved in a quadratic pair towards having Lie type. 
Its purpose may therefore be called \emph{group identification}.

As A. Premet and I. Suprunenko \cite{Premet-Suprunenko} put it, in \cite{Thompson} and \cite{Ho} ``groups generated by quadratic elements are classified as abstract finite groups and corresponding modules are not indicated explicitly.'' 
The article \cite{Premet-Suprunenko} by A. Premet and I. Suprunenko we just quoted attempts at remedying the lack of information on the module by listing finite groups of Lie type and representations thereof such that the pair they form is quadratic in J. Thompson's sense. This orthogonal line could conveniently be named \emph{representation zoology}.
Yet one then deals with a representation instead of a general module and this is much more accurate data.

As a matter of fact we shall adopt neither the \emph{group identification} nor the \emph{representation zoology} approach but a third one which qualifies as \emph{module linearization}: given a more-or-less concrete group of Lie type and an abstract module, can one retrieve a linear structure compatible with the action? Such a trend can be traced to a result of G. Glauberman \cite[Theorem 4.1]{Glauberman} which having among its assumptions both finiteness and quadraticity turns an abelian $p$-group into a sum of copies of the natural $\SL_2(\F_{p^n})$-module.
Following S. Smith \cite[Introduction]{Smith}, it was F.G. Timmesfeld who first asked whether similar results identifying the natural $\SL_2(\K)$-module among abstract quadratic modules would hold \emph{over possibly infinite fields}. As one sees this involves reconstructing a linear structure without the arsenal of finite group theory. Answers were given by F.G. Timmesfeld \cite[Proposition 2.7]{TGroups} and S. Smith \cite{Smith}.

Of course matters are a little more subtle than this rough historical account as one may be interested in \emph{simultaneous} identification: actually G. Glauberman \cite[Theorem 4.1]{Glauberman} also identified the group, and this combined direction has been explored extremely far by F.G. Timmesfeld \cite{Timmesfeld-Abstract-Quadratic}.

We shall follow the line of pure \emph{module linearization}. Our group or Lie ring is explicitly known to be $\SL_2(\K)$ or $\sl_2(\K)$; given a quadratic module, we wish to retrieve a compatible linear geometry. Parts of the present article, namely the Theme and Variations \ref{v:SL2quadratiqueCVi=CVG}--\ref{v:SL2quadratique}, are no original work but are adapted from F.G. Timmesfeld's book \cite{Timmesfeld}. What we add to the existing literature is the replacement of an assumption on the unipotent subgroup by an assumption on a single unipotent element, and the treatment of the Lie ring $\sl_2(\K)$.

The liminary digression ends here. Our main result is the following.
\medskip

\par\noindent\textbf{Variations \ref{v:d2=0quadratique}, \ref{v:SL2quadratique}, and \ref{v:sl2quadratique}.}
\textit{Let $\K$ be a field of characteristic $\neq 2, 3$, $\fG = \SL_2(\K)$ or $\sl_2(\K)$, and $V$ be a $\fG$-module. Suppose that there is a unipotent element $u$ (resp., nilpotent element $x$) of $\fG$ acting quadratically on $V$, meaning that $(u-1)^2$ or $x^2$ is zero in $\End V$. If\inmargin{Fixed} $\fG = \sl_2(\K)$ where $\K$ has characteristic $0$, suppose in addition that $V$ is $3$-torsion-free. Then $V$ is the direct sum of a $\fG$-trivial submodule and of copies of the natural representation $\fG$.}
\medskip

The result for the Lie ring $\sl_2(\K)$ (Variation \ref{v:sl2quadratique}) seems to be new. The result for the group $\SL_2(\K)$ is a non-trivial strengthening (Variation \ref{v:d2=0quadratique}) of F.G. Timmesfeld's work (Variation \ref{v:SL2quadratique}), as the assumption is now only about one unipotent \emph{element}, not about a unipotent subgroup; however the argument works only in characteristic $\neq 2, 3$.
\inmargin{Revised} It could be expected from J. Thompson's work in characteristic $\geq 5$ and Ho's delicate extension to characteristic $3$ (see the introductory digression above) that the case $p = 3$ would be quite harder if not different.

In the case of the Lie ring $\fG = \sl_2(\K)$, one can produce counter-examples in characteristic $3$\inmargin{Revised} but this requires the ``opposite'' nilpotent element $y$ to behave non-quadratically. In the case of the group $\fG = \SL_2(\K)$, I do not know.

The reader may also find of interest Variation \ref{v:algebrica}, whose lengthy proof indicates that reducing an $\SL_2(\K)$-module to an $\sl_2(\K)$-module is not any simpler than directly linearizing the former.
\medskip

The current section \S\ref{S:setting} is devoted to notations and basic observations. In \S\ref{S:quadratique} the core of the Smith-Timmesfeld argument for quadratic $\SL_2(\K)$-modules is reproduced; it will be generalized in following papers whence our present recasting it. Still on quadratic $\SL_2(\K)$-modules, \S\ref{s:centralizers} bears no novelty but \S\ref{s:length} may. In \S\ref{S:algebra}, the Lie ring $\sl_2(\K)$ and its quadratic modules are studied.
\medskip

\begin{notation*}
Let $\K$\inmargin{$\K$, $\fG$} be a field and $\fG$ be the $\K$-points of $\SL_2$ or $\sl_2$.
\end{notation*}

$\fG$ will thus denote either a group $G$ or a Lie ring $\fg$.

\subsection{The Group}

\begin{notation*}
Let $G$ be the group $\SL_2 (\K)$\inmargin{$G$}.
\end{notation*}

\begin{notation*}
For $\lambda \in \K$ (resp. $\K^\times$), let: \inmargin{$u_\lambda$, $t_\lambda$}
$$u_\lambda = \begin{pmatrix}1 & \lambda \\ 0 & 1\end{pmatrix} \quad \mbox{and}\quad t_\lambda = \begin{pmatrix}\lambda & 0 \\ 0 & \lambda^{-1}\end{pmatrix}$$
One simply writes $u = u_1$ and $i = t_{-1} \in Z(G)$\inmargin{$u$, $i$}.
\end{notation*}

If the characteristic is $2$, one has $i = 1$.

\begin{notation*}
Let\inmargin{$U, T$}:
$$U = \left\{u_\lambda: \lambda \in \K^+ \right\} \simeq \K^+
\quad\mbox{and}\quad
T = \left\{t_\lambda: \lambda \in \K^\times \right\} \simeq \K^\times$$
Let $B = U \rtimes T = N_G(U)$\inmargin{$B$}.
\end{notation*}

$B$ is a Borel subgroup of $G$ and $U$ is its unipotent radical, which is a maximal unipotent subgroup; $T$ is a maximal algebraic torus.

\begin{relations*}\
\begin{itemize}
\item
$u_\lambda u_\mu = u_{\lambda + \mu}$;
\item
$t_\lambda t_\mu = t_{\lambda\mu}$;
\item
$t_\mu u_\lambda t_{\mu^{-1}} = u_{\lambda \mu^2}$.
\end{itemize}
\end{relations*}

Note that in characteristic $\neq 2$, every element is a difference of two squares: consequently $\<T, u\> = T \ltimes U$.

\begin{notation*}
Let\inmargin{$w$} $w = \begin{pmatrix}0 & 1 \\ -1 & 0 \end{pmatrix}$.
\end{notation*}

\begin{relations*}
One has $w^2 = i$ and $w t_\lambda w^{-1} = t_{\lambda^{-1}} = t_{\lambda}^{-1}$.
\end{relations*}

\begin{relations*}
$u_\lambda wu_{\lambda^{-1}} wu_\lambda w= t_\lambda$,
and in particular $(uw)^3 = 1$.
\end{relations*}

The natural (left-) module $\Nat \SL_2(\K)$ corresponds to the natural action of $G$ on $\K^2$.

\subsection{The Lie Ring}

\begin{notation*}
Let $\fg$\inmargin{$\fg$} be the Lie ring $\sl_2(\K)$.
\end{notation*}

\begin{notation*}
For $\lambda \in \K$, let:\inmargin{$h_\lambda$, $x_\lambda$, $y_\lambda$}
$$h_\lambda = \left(
\begin{array}{cc}
 \lambda & 0 \\ 0 & -\lambda
\end{array}\right), \quad
x_\lambda = \left(
\begin{array}{cc}
 0 & \lambda \\ 0 & 0
\end{array}\right), \quad
y_\lambda = \left(
\begin{array}{cc}
 0 & 0 \\ \lambda & 0
\end{array}\right)$$
One simply writes $h = h_1$, $x = x_1$, $y = y_1$.\inmargin{$h, x, y$}
\end{notation*}

\begin{notation*}
Let\inmargin{$\fu, \ft$}:
$$\fu = \left\{x_\lambda: \lambda \in \K^+ \right\} \simeq \K^+
\quad\mbox{and}\quad
\ft = \left\{h_\lambda: \lambda \in \K^+ \right\} \simeq \K^+$$
Let $\fb = \fu \oplus \ft = N_\fg(\fu)$\inmargin{$\fb$}.
\end{notation*}

$\fb$ is a Borel subring of $\fg$ and $\fu$ is its nilpotent radical; $\ft$ is a Cartan subring.

\begin{relations*}\
\begin{itemize}
\item
$[h_\lambda, x_\mu] = 2 x_{\lambda \mu}$;
\item
$[h_\lambda, y_\nu] = - 2 y_{\lambda \nu}$;
\item
$[x_\mu, y_\nu] = h_{\mu\nu}$.
\end{itemize}
\end{relations*}

The natural (left-) module $\Nat \sl_2(\K)$ corresponds to the natural action of $\fg$ on $\K^2$.

\subsection{The Module}

\begin{notation*}
Let $V$\inmargin{$V$} be a $\fG$-module, that is a $G$- or $\fg$-module.
\end{notation*}

The names of the elements of $\fG$ will still denote their images in $\End V$.

\begin{notation*}
When $\fG = G$, one lets for $\lambda \in \K$: \inmargin{$\d_\lambda$, $\d$}$\d_\lambda = u_\lambda - 1 \in \End V$. One simply writes $\d = \d_1$.
\end{notation*}

\begin{relations*}\
\begin{itemize}
\item
$\d_\lambda \circ \d_\mu = \d_\mu \circ \d_\lambda$;
\item
$t_\lambda \d_\mu = \d_{\lambda^2\mu} t_\lambda$;
\item
$\d_{\lambda+\mu} = \d_\lambda + \d_\mu + \d_\lambda \circ \d_\mu$.
\end{itemize}
\end{relations*}
\begin{proofclaim}
The first claim is by abelianity of $U$; the second comes from the action of $T$ on $U$. Finally, denoting by $u_\lambda$ the corresponding element in the group ring (or more precisely its image in $\End V$), one has:
\begin{eqnarray*}
\d_{\lambda + \mu} & = & u_{\lambda + \mu} -1 = u_\lambda u_\mu - 1 = (u_\lambda u_\mu - u_\mu) + (u_\mu - 1)\\
& = & \d_\lambda u_\mu + \d_\mu = \d_{\lambda} (\d_\mu + 1) + \d_\mu = \d_\lambda \d_\mu + \d_\lambda + \d_\mu
\end{eqnarray*}
as desired.
\end{proofclaim}

\begin{notation*}
When $\fG = \fg$, one lets for $i \in \Z$: $E_i(V) = \{a \in V: h\cdot a = iv\}$\inmargin{$E_i(V)$}. When there is no ambiguity on the module, one simply writes $E_i$.
\end{notation*}

Each $h_\lambda$ (resp. $x_\mu$, resp. $y_\nu$) maps $E_i$ into $E_i$ (resp. $E_{i+2}$, resp. $E_{i-2}$).
\inmargin{Added} One should however be careful that if the module contains torsion, the various $E_i$'s need not be in direct sum.

\begin{notation*}
The length \inmargin{$\ell_U(V), \ell_\fu(V)$} of $V$ is the smallest integer, if there is one:
\begin{itemize}
\item
when $\fG = G$, such that $[U, \dots, U, V] = 0$ ($U$-length);
\item
when $\fG = \fg$, such that $\fu \dots \fu \cdot V = 0$ ($\fu$-length).
\end{itemize}
A length $2$ module is called quadratic.
\end{notation*}

Clearly, if $V$ is simple (i.e. without a proper, non-trivial $\fG$-submodule), then $V$ either has prime exponent, or is torsion-free and divisible. We shall not always assume this.

The group $G$ is said to act trivially on $V$ if it centralizes it, that is if the image of $G$ in $\End V$ is $\{\Id\}$; the Lie ring $\fg$ is said to act trivially on $V$ if it annihilates it, that is if the image of $\fg$ in $\End V$ is $\{0\}$.
We then say that $V$ is $G$- (respectively $\fg$-) trivial.
The following observations will be used with no reference.

\begin{observation*}
Suppose that $\fG = \fg = \sl_2(\K)$. Let $V$ be a $\fg$-module.
\begin{enumerate}
\item
If $\K$ has characteristic $p$ and $V$ is $p$-torsion-free, then $V$ is $\fg$-trivial.
\item
If $\K$ has characteristic $0$ and $V$ is torsion, then $V$ is $\fg$-trivial.
\end{enumerate}
\end{observation*}
\begin{proofclaim}
Fix $a \in V\setminus \{0\}$ and any element $z$ of $\fg$.
\begin{enumerate}
\item
If $\K$ has characteristic $p$, then $\fg$ has exponent $p$. Suppose that $V$ is $p$-torsion-free; then $pz \cdot a = 0$ implies that $z \cdot a = 0$: $\fg$ annihilates $V$.
\item
If $\K$ has characteristic $0$, then $\fg$ is divisible. Suppose that $V$ is torsion; let $n$ be the order of $a$. Then $n \left(\frac{1}{n}z\right) \cdot a = 0 = z \cdot a$: $\fg$ annihilates $V$.\qedhere
\end{enumerate}
\end{proofclaim}

The case of the group is hardly less trivial.

\begin{observation*}
Suppose that $\fG = G = \SL_2(\K)$. Let $V$ be a $G$-module. Suppose that $\d$ is nilpotent in $\End V$.
\begin{enumerate}
\item
If $\K$ has characteristic $p$ and $V$ is $p$-torsion-free, then $V$ is $G$-trivial.
\item
If $\K$ has characteristic $0$ and $V$ is torsion, then $V$ is $G$-trivial.
\end{enumerate}
\end{observation*}
\begin{proofclaim}\
\begin{enumerate}
\item
We show that $u$ centralizes $V$.
Otherwise, from the assumptions, there is $a_2 \in \ker \d^2 \setminus \ker \d$. Let $a_1 = \d(a_2) \in \ker \d \setminus \{0\}$.
Since $\K$ has characteristic $p$ and $V$ is $p$-torsion-free, $u^p \cdot a_2 = a_2 = a_2 + p a_1$ implies $a_1 = 0$: a contradiction. Hence $u$ centralizes $V$.
So $C_G(V)$ is a normal subgroup of $G$ containing an element of order $p$: it follows that $C_G(V) = G$ (this still holds of $\K = \F_2$ or $\F_3$).
\item
If $\K$ has characteristic $0$, the previous argument is no longer valid.
Since $V$ splits as the direct sum of its $p$-torsion components, one may assume that $V$ is a $p$-group. We further assume that $V$ has exponent $p$.

For any $a \in V$ and any integer $k$, one has $u^k \cdot a = \sum_{i \leq k} \binom{k}{i} \d^i (a)$; since $\d$ has finite order $\ell$, for $k \geq \ell$ one even has $u^k \cdot a = \sum_{i = 0}^\ell \binom{k}{i} \d^i (a)$. But since $V$ has exponent $p$, for $k$ big enough (independently of $a$) one finds $u^k \cdot a = a$.

Hence $u^k$ centralizes $V$.
Here again the normal closure of $u^k$ is $G$, which must centralize $V$. We finish the argument. Let $V_{p^n}$ be the $G$-submodule of $V$ of exponent $p^n$. Then $G$ centralizes every $V_{p^n}/V_{p^{n-1}}$. But $G = \SL_2(\K)$ is perfect; it therefore centralizes $V$.
\qedhere
\end{enumerate}
\end{proofclaim}

\section{The Natural Module}\label{S:quadratique}

\begin{thema*}\label{v:thema}
Let $\K$ be a field, $G = \SL_2(\K)$, and $V$ be a simple $G$-module of $U$-length $2$. Then there exists a $\K$-vector space structure on $V$ making it isomorphic to $\Nat \SL_2(\K)$.
\end{thema*}

This theorem was proved by F.G. Timmesfeld in a more general context (Theorem 3.4 of chapter I in his book \cite{Timmesfeld}) and independently by S. Smith \cite{Smith}. Let us adapt the proof to our notations.

\begin{proof}
The assumption means that $[U, U, V] = 0$. Let $Z_1 = C_V(U)$, so that $U$ centralizes $V/Z_1$. Recall that one lets $\d_\lambda = u_\lambda - 1 \in \End V$. These functions map $V$ to $Z_1$ and annihilate $Z_1$.

Observe that by simplicity, $C_V(G) = 0$.

\subsection{Finding a Decomposition}

\begin{step}
$Z_1 \cap w \cdot Z_1 = 0$.
\end{step}
\begin{proofclaim}
$Z_1 \cap w \cdot Z_1 = C_V(U, w U w^{-1}) = C_V (G) = 0$.
\end{proofclaim}

Recall that $i$ denotes the central element of $G$ ($i = 1$ in characteristic $2$).

\begin{step}
For all $a_1 \in Z_1$, $\d_\lambda (w \cdot a_1) = i t_\lambda \cdot a_1$.
\end{step}
\begin{proofclaim}
Let $b_1 = \d_\lambda (w \cdot a_1)$ and $c_1 = \d_{\lambda^{-1}} (w\cdot b_1)$; by assumption, $b_1$ and $c_1$ lie in $Z_1$. Then:
\begin{eqnarray*}
u_{\lambda^{-1}} w u_\lambda w \cdot a_1 & = & u_{\lambda^{-1}} w\cdot (w\cdot a_1 + b_1) \\
& = & i \cdot a_1 + w \cdot b_1 + c_1\\
 = (u_\lambda w)^{-1} t_\lambda \cdot a_1 & = & w^{-1} u_{-\lambda} t_\lambda \cdot a_1 \\
& = & i w t_\lambda \cdot a_1
\end{eqnarray*}
So $i \cdot a_1 + c_1 = w \cdot (i t_\lambda \cdot a_1 - b_1) \in Z_1\cap w \cdot Z_1 = 0$ and the claim follows.
\end{proofclaim}

\subsection{Linear Structure}

\begin{localnotation}
For $\lambda \in \K$ and $a_1 \in Z_1$, let:
$$\left\{ \begin{array}{rcl}
      \lambda \cdot a_1 & = & t_\lambda \cdot a_1 \\
\lambda \cdot (w\cdot a_1) & = & w\cdot(\lambda \cdot a_1)
     \end{array}
\right.$$
\end{localnotation}

\begin{step}
This defines an action of $\K$ on $Z_1\oplus w \cdot Z_1$.
\end{step}
\begin{proofclaim}
It is clearly well-defined. The action is obviously multiplicative on $Z_1 \oplus w \cdot Z_1$, because each term is $T$-invariant. Moreover one has:
\begin{itemize}
\item
on $Z_1$:
\begin{eqnarray*}
(\lambda + \mu) \cdot a_1 & = & t_{\lambda + \mu} \cdot a_1 = i \cdot \d_{\lambda + \mu} (w \cdot a_1)\\
& = & i \cdot (\d_\lambda (w \cdot a_1) + \d_\mu (w \cdot a_1) + \d_\lambda \d_\mu (w \cdot a_1))\\
& = & i \cdot \d_\lambda (w \cdot a_1) + i \cdot \d_\mu (w \cdot a_1)
= t_\lambda \cdot a_1 + t_\mu \cdot a_1\\
& = & \lambda \cdot a_1 + \mu \cdot a_1
\end{eqnarray*}
\item
on $w \cdot Z_1$:
\begin{eqnarray*}
(\lambda + \mu) \cdot (w \cdot a_1) & = & w \cdot ((\lambda + \mu) \cdot a_1) = w \cdot (\lambda \cdot a_1 + \mu \cdot a_1)\\
& = & w \cdot (\lambda \cdot a_1) + w \cdot (\mu \cdot a_1)\\
& = & \lambda \cdot (w \cdot a_1) + \mu \cdot (w \cdot a_1)
\end{eqnarray*}
\end{itemize}
and everything is proved.
\end{proofclaim}

\begin{step}
$G$ is linear on $Z_1 \oplus w \cdot Z_1$.
\end{step}
\begin{proofclaim}
Clearly $\<T , w\>$ acts linearly. Moreover $\d_\lambda$ is trivially linear on $Z_1$. Finally $\d_\lambda (\mu \cdot (w\cdot a_1)) = \d_\lambda (w\cdot(\mu \cdot a_1)) = i t_\lambda\cdot (\mu \cdot a_1) = \mu \cdot (i t_\lambda \cdot a_1) = \mu \cdot \d_\lambda (w\cdot a_1)$
so $\d_\lambda$ is linear, and $u_\lambda$ is therefore too.
\end{proofclaim}

$V$ being simple is additively generated by the $G$-orbit of any $a_1 \in Z_1 \setminus \{0\}$, and one then sees that $V \simeq \K^2$ as the natural $G$-module.
This finishes the proof.
\end{proof}

\begin{remark*}
Note that although there are a priori several $\K$-vector spaces structures such that $G$ acts linearly (twist the action by any field automorphism), our construction is uniquely defined. It is functorial: if $V_1$ and $V_2$ are two simple $\SL_2(\K)$-modules and $\varphi: V_1 \to V_2$ is a morphism of $\SL_2(\K)$-modules, then for our construction $\varphi$ is $\K$-linear.
\end{remark*}

\section{First Variations}

\subsection{Centralizers}\label{s:centralizers}

The statements of this subsection can be found in F.G. Timmesfeld's book \cite{Timmesfeld}.

\begin{variation}\label{v:SL2quadratiqueCVi=CVG}
Let $\K$ be a field of characteristic $\neq 2$ with more than three elements, $G = \SL_2 (\K)$, and $V$ be a $G$-module. Suppose that $V$ has $U$-length at most $2$. Then $G$ centralizes $C_V(i)$.
\end{variation}
\begin{proof}
We assume that the central involution $i$ centralizes $V$ and show that $G$ does too. By our assumptions, $G$ centralizes the $2$-torsion component $V_2$ of $V$. Recall that one writes $\d = \d_1$.

Let $a \in C_V(U)$, $b = \d (w \cdot a)$, and $c = \d (w \cdot b)$. One then has:
\begin{eqnarray*}
u w u w \cdot a & = & uw \cdot (w \cdot a + b)\\
& = & a + w \cdot b + c\\
= w^{-1} u \cdot a & = & w \cdot a
\end{eqnarray*}
Hence $w\cdot (a - b) = a + c$. But by assumption on the $U$-length, $b, c \in C_V(U)$, so $b-a \in C_V(U, wUw^{-1}) = C_V(G)$. Let us resume:
\begin{eqnarray*}
u w u w \cdot a & = & uw \cdot (w \cdot a + (b-a) + a)\\
& = & a + (b-a) + (w \cdot a + (b-a) + a)\\
= w^{-1} u \cdot a & = & w \cdot a\end{eqnarray*}
Therefore $2 b = 0$, that is $b \in V_2 \leq C_V(G)$, whence $a = (a-b) + b \in C_V(G)$.

As a conclusion $G$ centralizes $C_V(U)$. But by assumption on the $U$-length, $U$ centralizes $V/C_V(U)$, so by the same argument $G$ centralizes $V/C_V(U)$ as well. Now $G$ being perfect by the assumptions on $\K$, $G$ does centralize $V$.
\end{proof}

\begin{variation}[{\cite[Lemma 3.1 of chapter I]{Timmesfeld}}]\label{v:SL2quadratiquecoherente}
Let $\K$ be a field of characteristic $\neq 2$ having more that three elements, $G = \SL_2 (\K)$, and $V$ be a $G$-module of $U$-length $2$ satisfying $C_V(G) = 0$. Then for any $\lambda \in \K^\times$, $[u_\lambda, V] = [U, V] = C_V(U) = C_V(u_\lambda)$. In particular $C_V(u_\lambda)$ does not depend on $\lambda$.
\end{variation}
\begin{proof}
We can prove it as a Corollary to the Theme (modulo a few adjustements) or argue as follows. Since $C_V(G) = 0$, $V$ is $2$-torsion-free, and by Variation \ref{v:SL2quadratiqueCVi=CVG}, $C_V(i)= 0$. It follows that $i$ inverts $V$.

By assumption on the $U$-length, $[u_\lambda, V] \leq [U, V] \leq C_V(U) \leq C_V(u_\lambda)$. Let $a \in C_V(u_\lambda)$: we show that $a \in [u_\lambda, V]$.
Let $b = \d_{\lambda^{-1}} (w \cdot a)$ and $c = \d_\lambda (w \cdot b)$, so that:
\begin{eqnarray*}
u_\lambda w u_{\lambda^{-1}}w \cdot a & = & (u_{\lambda} w) \cdot (w \cdot a + b)\\
& = & - a + w \cdot b + c\\
= (w^{-1} u_{\lambda^{-1}}^{-1} t_{\lambda^{-1}}) \cdot a
& = & - (w t_{\lambda^{-1}} u_{-\lambda}) \cdot a
\end{eqnarray*}
Hence $a - c = w \cdot (b + t_{\lambda^{-1}} \cdot a)$. But on the one hand $c \in [U, V] \leq C_V(u_\lambda)$, so $a-c$ commutes with $u_\lambda$, and on the other hand $t_{\lambda^{-1}} \cdot a \in C_V(t_{\lambda^{-1}} u_\lambda t_\lambda) = C_V(u_{\lambda^{-1}})$ so $a-c$ also commutes with $w u_{\lambda^{-1}} w^{-1}$. Hence $C_G(a-c)$ contains:
\begin{eqnarray*}
(u_\lambda w u_{\lambda^{-1}} w^{-1})^3 & = & i (u_\lambda w u_{\lambda^{-1}} w)^3\\
& = & i (u_\lambda w u_{\lambda^{-1}} w u_\lambda w) (u_{\lambda^{-1}} w u_\lambda w u_{\lambda^{-1}} w)\\
& = & i t_\lambda t_{\lambda^{-1}}
\end{eqnarray*}
So $i$ which inverts $V$, centralizes $a - c$; since $V$ is $2$-torsion-free it follows that $a = c \in [u_\lambda, V]$.
\end{proof}

Recall that when $i$ is an involutive automorphism of an abelian group $V$, one lets $V^{+_i} = \{v \in V: i \cdot v = v\}$ and $V^{-_i} = \{v \in V: i \cdot v = -v\}$; when there is no ambiguity one simply writes $V^+$ and $V^-$. If $V$ is $2$-torsion-free then $V^+ \cap V^- = 0$; if $V$ is $2$-divisible then $V = V^+ + V^-$. Actually if $[i, V]$ is $2$-divisible, one has $V = V^+ + [i, V]$.

\begin{variation}[{\cite[Exercise 3.8.1 of chapter I]{Timmesfeld}}]\label{v:SL2quadratique}
Let $\K$ be a field of characteristic $\neq 2$ with more than three elements, $G = \SL_2 (\K)$, and $V$ be a $G$-module of $U$-length $\leq 2$. Then $V = C_V(G) \oplus [G, V]$, and there exists a $\K$-vector space structure on $[G, V]$ making it isomorphic to a direct sum of copies of $\Nat \SL_2(\K)$. In particular $C_V(U) = C_V(u_\lambda)$ for any $\lambda \in \K^\times$.
\end{variation}
\begin{proof}
We have made no assumption on $2$-divisibility or $2$-torsion-freeness of $V$, so one may not a priori decompose $V$ as $V^+$ and $V^-$ under the action of the central involution; the argument is more subtle.

By Variation \ref{v:SL2quadratiqueCVi=CVG}, $G$ centralizes $V^+$, that is $V^+ = C_V(G)$. Let $W = [G, V]$ and $\overline{W} = W/C_W(G) = W/W^+$; these are $G$-modules of $U$-length $\leq 2$.
By perfectness of $G$, $C_{\overline{W}}(G) = 0$.

One then reads the proof of the Theme again and sees that simplicity was only used to show that $C_V(G) = 0$.
In particular the Theme constructs, for any $\bar{a}_1 \in C_{\overline{W}}(U) \setminus \{0\}$, a $\K$-linear structure on $\<G \cdot \bar{a}_1\>$ such that $G$ acts naturally.
We then take a maximal family of such vector planes in direct sum.
By perfectness of $G$ one has $W = [G, W]$ and $\overline{W} = [G, \overline{W}]$.
Since $G = \<U, w U w^{-1}\>$, one has $\overline{W} = [G, \overline{W}] = [U, \overline{W}] + [w U w^{-1}, \overline{W}] \leq \<G \cdot C_{\overline{W}}(U)\>$, so $\overline{W}$ is itself a direct sum of vector planes all isomorphic to the natural representation of $G$.

In particular $i$ inverts $\overline{W}$, and the characteristic of $\K$ being $\neq 2$, $\overline{W}$ is $2$-divisible and $2$-torsion-free.
Let $a \in W$. As $\overline{W}$ is $2$-divisible, there is $b \in W$ such that $a - 2b \in C_W(G)$. Since $i$ inverts $\overline{W}$, $(i+1) \cdot b \in C_W(G)$. We take the sum: $a + (i-1) \cdot b \in C_W(G)$. This means that $W \leq [i, W] + C_W(G)$, and therefore $W = [G, W] \leq [G, [i, W]] = [i, W]$.

Now let $a \in C_W(G) = W^+$; as $W = [i, W]$ there is $b \in W$ such that $a = i\cdot b - b$, and applying $i$ one gets $2b \in C_W(G)$. But $\overline{W}$ is $2$-torsion-free, so $b \in C_W(G)$ and $a = 0$.
This implies $C_W(G) = 0$, and retrospectively $\overline{W} = W = [i, W] = [i, V] = [G, V]$ which is $2$-divisible and $2$-torsion-free. One thus has $V = V^+ + [i, V] = C_V(G) \oplus [G, V]$.

The final claim on centralizers is obtained by Variation \ref{v:SL2quadratiquecoherente}, or more prosaically by inspection in each copy of $\Nat \SL_2(\K)$.
\end{proof}

\subsection{Length}\label{s:length}

\begin{variation}\label{v:dk=0dlambda2k-1=0}
Let $\K$ be a field of characteristic $\neq 2$, $G = \SL_2(\K)$, and $B$ be a Borel subgroup of $G$. Let $V$ be a $B$-module. Suppose that $V$ has $u$-length at most $k$, meaning that $\d^k = 0$ in $\End V$. Then for any $\lambda \in \K$, $\d_\lambda^{2k-1} = 0$.
\end{variation}
\begin{proof}
Indeed, $\lambda$ is a difference of two squares $\lambda = \mu^2 - \nu^2$, so $\d_\lambda = \d_{\mu^2 - \nu^2} = \d_{\mu^2} + \d_{-\nu^2} + \d_{\mu^2}\d_{-\nu^2}$.
But $\d_{\mu^2}$ and $\d_{\nu^2}$ are $T$-conjugate to $\d$, so they have order at most $k$. Moreover $\d_{-\nu^2} = - \d_{\nu^2} + \d_{\nu^2}^2 + \dots + (-1)^{k-1} \d_{\nu^2}^{k-1}$. It is now clear that $\d_\lambda^{2k-1} = 0$.
\end{proof}


\begin{variation}\label{v:repartitioncarres}
Let $\K$ be a field, $G = \SL_2 (\K)$, $B \leq G$ be a Borel subgroup, and $V$ be a $B$-module. Suppose that $V$ has $u$-length $\leq k$, meaning $\d^k = 0$. If every element of $\K$ is a (positive or negative) integer multiple of a square, then for every $\lambda \in \K$, one has $\d_\lambda^k = 0$.
\end{variation}
\begin{proof}
Let $\lambda$ be a square. Then $\d_\lambda$ is $T$-conjugate to $\d$, so $\d_\lambda^k = 0$. Now for any $n \in \N$, $\d_{n \lambda} = u_{n \lambda} - 1 = u_\lambda^n - 1 = \sum_{j = 1}^n \binom{n}{j} \d_\lambda^j$\inmargin{Fixed typo}, the $k^{\rm th}$ power of which is zero.
Finally $u_{-_\lambda} = u_\lambda^{-1} = (1 + \d_\lambda)^{-1} = 1 - \d_\lambda + \d_\lambda^2 \dots + (-1)^{k-1} \d^{k-1}$, so $\d_{-\lambda}^k = 0$. Hence any integer multiple of $\lambda$ will satisfy $\d_\lambda^k = 0$. Our assumption is precisely that every element of $\K$ is of this form.
\end{proof}

\begin{variation}\label{v:dlambdan=0longueurn}
Let $\K$ be a field, $G = \SL_2(\K)$, and $U$ be a unipotent subgroup of $G$. Let $V$ be a $U$-module. If for all $\lambda \in \K$, $\d_\lambda^n = 0$ in $\End V$ and $V$ is $n!$-torsion-free, then $V$ has $U$-length $\leq n$.
\end{variation}
\begin{proof}
Suppose that for any $\lambda$, one has $\d_\lambda^n = 0$; we show that every product $\d_{\lambda_1} \dots \d_{\lambda_n}$ annihilates $V$.
Fix $\lambda$ and $\mu$. Then $\d_{\lambda + \mu} = \d_\lambda + \d_\mu + \d_\lambda\d_\mu$ and $\d_{\lambda+\mu}^n = 0$, so that:
$$0 = \sum_{i = 0}^n \binom{n}{i} (\d_\lambda\d_\mu)^{n-i} \sum_{j = 0}^i \binom{i}{j} \d_\lambda^j \d_\mu^{i-j}
= \sum_{0 \leq j \leq i \leq n} \binom{n}{i} \binom{i}{j} \d_\lambda^{n-i+j} \d_\mu^{n-j}
$$
The monomials occuring in this sum have weight $2n-i$. We show by induction on $k = 2n-1 \dots n$ that every monomial of weight $\geq k$ is zero. When $k = 2n-1$, the only two such monomials are $\d_\lambda^n \d_\mu^{n-1}$ and $\d_\lambda^{n-1} \d_\mu^n$: both are zero by assumption.

So suppose the result holds for $k+1$; we prove it for $k$, with $k \geq n$. Multiplying the equation by $\d_\lambda^{k-n}$, one finds:
$$0 = \sum_{0 \leq j \leq i \leq n} \binom{n}{i} \binom{i}{j} \d_\lambda^{k-i+j} \d_\mu^{n-j}$$
But when $i < n$, the terms have weight $n+k-i \geq k +1$, so all monomials are zero. Hence only the terms with $i = n$ remain, that is:
\begin{equation*}\label{e:dlambdak=0Uk=0}
0 = \sum_{j = 0}^n \binom{n}{j} \d_\lambda^{k-n+j} \d_\mu^{n-j} = \sum_{j = 1}^{n-1} \binom{n}{j} \d_\lambda^{k-n+j} \d_\mu^{n-j}
\end{equation*}

We now replace $\mu$ by $i\mu$. Since $\d_{i \mu}$ is equal to $i \d_\mu$ modulo terms of weight $\geq 2$, one actually has for all $i = 1 \dots n-1$:
$$0 = \sum_{j = 1}^{n-1} \binom{n}{j} i^{n-j} \d_\lambda^{k-n+j} \d_\mu^{n-j}$$

This gives $n-1$ equations in $n-1$ variables, with determinant:
$$\left|\binom{n}{j} i^{n-j}\right|_{i, j = 1 \dots n-1} = \left|i^{j}\right|_{i, j = 1 \dots n-1} \prod_{j=1}^{n-1} \binom{n}{j} = \prod_{j=1}^{n-1} j! \binom{n}{j} = \frac{(n!)^{n-1}}{\prod_{j = 1}^{n-1} (n-j)!}$$
Since $V$ is $n!$-torsion-free, one deduces that all terms are trivial: the latter are the monomials of weight $k$.

This completes the induction. It follows in particular that $\d_\lambda^{n-1} \d_\mu$ is trivial in $\End V$. But $\mu$ being fixed, $U$ acts on $\im \d_\mu$ which is $(n-1)!$-torsion-free, and $\d_\lambda^{n-1}$ acts trivially. By induction on $n$, one gets that every product $\d_{\mu_n} \dots \d_{\mu_1}$ is trivial on $V$: which was to be proved.
\end{proof}

\begin{remark*}
If $\K$ has characteristic $p$ and $V$ has exponent $p$, then without any assumption on $u$, every unipotent element $u_\lambda$ acts with length at most $p$: one has indeed $u_\lambda^p = 1 = (1 + \d_\lambda)^p = 1 + \d_\lambda^p$. Yet $V$ does not necessarily have $U$-length at most $p$, even if $V$ actually is a $G$-module.

For any prime $p$, one may check that the Steinberg module $\St \SL_2(\F_{p^2})$ is a simple $\SL_2(\F_{p^2})$-module of exponent $p$ with $U$-length $> p$: all unipotent elements have length $p$, but the action hasn't. Going to $\St \SL_2(\F_{p^n})$ one can even make the $U$-length arbitrarily big.
\end{remark*}

\begin{variation}\label{v:d2=0quadratique}
Let $\K$ be a field of characteristic $\neq 2, 3$, $G = \SL_2(\K)$, and $V$ be a $G$-module. Suppose that $V$ has $u$-length $\leq 2$, meaning that $\d^2 = 0$. Then $V$ has $U$-length $\leq 2$.
\end{variation}
\begin{proof}
\begin{step}\label{st:d2=0quadratique:CVG=0}
We may assume $C_V(G) = 0$.
\end{step}
\begin{proofclaim}
Let $\bar{V} = V/C_V(G)$; by perfectness of $G$, $C_{\bar{V}}(G) = 0$, and one still has $\d^2 = 0$ in $\End \bar{V}$. Suppose the result is proved for $\bar{V}$; we shall prove it for $V$.

Since $\bar{V}$ has $U$-length at most $2$ and $C_{\bar{V}}(G) = 0$, $\bar{V}$ is by Variation \ref{v:SL2quadratique} a direct sum of copies of the natural representation of $G$. In particular the central involution $i$ inverts $\bar{V}$, which is $2$-divisible and $2$-torsion-free.
For any $a \in V$, there is therefore $b \in V$ such that $a - 2b \in C_V(G)$; moreover $(1 + i) \cdot b \in C_V(G)$, so $a + (i-1) \cdot b \in C_V(G)$, proving $V = [i, V] + C_V(G)$.
Now let $a\in [i, V] \cap C_V(G)$. Then there is $b \in V$ such that $a = [i, b] \in C_V(G) \leq C_V(i)$, so $2b \in C_V(i)$. Since $i$ inverts $\bar{V}$, $(i+1) \cdot 2b = 4 b \in C_V(G)$, and as $\bar{V}$ is $2$-torsion-free, $b \in C_V(G) \leq C_V(i)$, whence $a = 0$.

One thus has $V = [i, V] \oplus C_V(G)$. In particular $[i, V] \simeq \bar{V}$ as $G$-modules, and $V$ has $U$-length $\leq 2$.
\end{proofclaim}

It follows from the assumptions on the base field that $V$ is $6$-torsion-free.
By Variations \ref{v:dk=0dlambda2k-1=0} and \ref{v:dlambdan=0longueurn}, $V$ has $U$-length at most $3$: $[U, U, U, V] = 0$. Let $Z_1 = C_V(U)$ and $Z_2$ be defined by $Z_2/Z_1 = C_{V/Z_1} (U)$. These subgroups are $B$-invariant; the $\d_\lambda$'s map $V$ into $Z_2$, $Z_2$ into $Z_1$, and annihilate $Z_1$. We must show that $Z_2 = V$.

\begin{step}
$C_V(i) = 0$.
\end{step}
\begin{proofclaim}
Consider $C_V(i)$ which is $G$-invariant and satisfies our assumptions; we may therefore suppose $V = C_V(i)$.
Let $a_1 \in Z_1$, $b_2 = \d (w \cdot a_1)$, and $c_2 = \d (w \cdot b_2)$. Note that $b_2, c_2 \in \im \d \leq \ker \d$. Then:
\begin{eqnarray*}
(uw)^{-1} \cdot a_1 & = & uwuw \cdot a_1\\
= w u^{-1} \cdot a_1 & = & uw \cdot (w\cdot a_1 + b_2)\\
= w \cdot a_1 & = & a_1 + w \cdot b_2 + c_2
\end{eqnarray*}
We apply $\d$: since $c_2 \in \im \d \leq \ker \d$, there remains $b_2 = c_2$. In particular $(w-1) \cdot a_1 = (w+1)\cdot b_2$. We apply $(w-1)$: one finds $(w-1)^2 \cdot a_1 = 2 (1-w) \cdot a_1 = 0$. Since $V$ is $2$-torsion-free, one has $w \cdot a_1 = a_1 \in Z_1 \cap w \cdot Z_1 = C_V(U, w U w^{-1}) = C_V(G) = 0$.
Hence $Z_1 = 0$, and since $V$ has finite $U$-length, $V = 0$.
\end{proofclaim}

In particular (and with no assumptions on $2$-divisibility), $i$ inverts $V$.

\begin{localnotation}
For any $\lambda \in \K^\times$, let $f_\lambda: Z_2 \to Z_2$ be such that $f_\lambda (a_2) = \d_\lambda (w\cdot a_2)$.
\end{localnotation}

It is not clear a priori whether $f_\lambda$ stabilizes $Z_1$.

\begin{step}\label{st:d2=0quadratique:naturelles}
If $a_1 \in Z_1 \cap w \cdot Z_2$ and $\lambda \in \K^\times$, then $f_\lambda (a_1) = - t_\lambda \cdot a_1$.
\end{step}
\begin{proofclaim}
For any $g \in G$,
\begin{itemize}
\item
either $g \in B$, in which case $g \cdot a_1 \in Z_1 \leq Z_2$;
\item
or $g \in B w U$, in which case $g\cdot a_1 \in Z_2$.
\end{itemize}
Let $V_0 = \<G\cdot (Z_1 \cap w \cdot Z_2)\>$: $V_0$ is therefore a $G$-submodule of $V$ included in $Z_2$, whence of $U$-length $\leq 2$. By Variation \ref{v:SL2quadratique} and since the involution inverts $V$, $V_0$ is a direct sum of copies of the natural representation of $G$.
It follows that for all $a_1 \in Z_1 \cap w \cdot Z_2$, $f_\lambda (a_1) = \d_\lambda (w \cdot a_1) = - t_\lambda \cdot a_1$.
\end{proofclaim}

We now go to the group ring $\Z[G]$, or more precisely its image in $\End (V)$. We shall drop parentheses and the application point $\cdot$ of a function to an element. There is no risk of confusion.

\begin{step}\label{st:d2=0quadratique:formules}
For any $\mu \in \K^\times$ and $a_2 \in Z_2$, one has:
\begin{eqnarray}
t_\mu f_{\mu^{-1}} a_2 & = & f_\mu t_{\mu^{-1}} a_2
\label{eq:d2=0quadratique:formule1}\\
- a_2 - \d_{\mu^{-1}} a_2 + w f_\mu a_2 + f_{\mu^{-1}} f_\mu a_2 & = & - w t_\mu a_2 + w t_\mu \d_{\mu^{-1}} a_2
\label{eq:d2=0quadratique:formule2}\\
- \d_{\mu^{-1}} a_2 + f_{\mu^{-1}} f_\mu a_2 + \d_{\mu^{-1}} f_{\mu^{-1}} f_\mu a_2 & = & - t_{\mu^{-1}} f_\mu a_2 + t_{\mu^{-1}} f_\mu \d_{\mu^{-1}} a_2
\label{eq:d2=0quadratique:formule3}
\end{eqnarray}
\end{step}
\begin{proofclaim}
First of all:
$$t_\mu f_{\mu^{-1}} a_2 = t_\mu \d_{\mu^{-1}} w a_2 = \d_\mu t_\mu w a_2 = \d_\mu w t_{\mu^{-1}} a_2 = f_\mu t_{\mu^{-1}} a_2$$

This proves (\ref{eq:d2=0quadratique:formule1}), which we shall use with no reference. Now to (\ref{eq:d2=0quadratique:formule2}).
On the one hand $u_\mu w a_2 = w a_2 + f_\mu a_2$, and since $a_2 \in Z_2$, one has on the other hand $u_{-\mu^{-1}} a_2 = a_2 - \d_{\mu^{-1}} a_2$, so that:
\begin{eqnarray*}
u_{\mu^{-1}} w u_\mu w a_2 & = & u_{\mu^{-1}} w (w a_2 + f_\mu a_2) = - u_{\mu^{-1}} a_2 + w f_\mu a_2 + f_{\mu^{-1}} f_\mu a_2\\
& = & - a_2 - \d_{\mu^{-1}} a_2 + w f_\mu a_2 + f_{\mu^{-1}} f_\mu a_2\\
= (u_\mu w)^{-1} t_\mu a_2 & = & - w u_{-\mu} t_\mu a_2 = - w t_\mu u_{-\mu^{-1}} a_2 = - w t_\mu a_2 + w t_\mu \d_{\mu^{-1}} a_2
\end{eqnarray*}
which proves (\ref{eq:d2=0quadratique:formule2}). To derive (\ref{eq:d2=0quadratique:formule3}), apply $\d_{\mu^{-1}}$.
\end{proofclaim}

\begin{step}\label{st:d2=0quadratique:court-circuit}
If $b_2 \in Z_2$ and $\lambda \in \K^\times$ are such that:
$$\left\{ \begin{array}{l}
f_{\lambda^{-1}} f_\lambda b_2 = - t_{\lambda^{-1}} f_\lambda b_2 + t_{\lambda^{-1}} \d_\lambda f_\lambda b_2\\
\d_\lambda f_\lambda b_2 \in w \cdot Z_2
\end{array}\right.$$
then $\d_\lambda f_\lambda b_2 = 0$.
\end{step}
\begin{proofclaim}
We apply formula (\ref{eq:d2=0quadratique:formule2}) of Step \ref{st:d2=0quadratique:formules} with $a_2 = f_\lambda (b_2)$ and $\mu = \lambda^{-1}$:
$$ - f_\lambda b_2 - \d_\lambda f_\lambda b_2 + w f_{\lambda^{-1}} f_\lambda b_2 + f_\lambda f_{\lambda^{-1}} f_\lambda b_2 = - w t_{\lambda^{-1}} f_\lambda b_2 + w t_{\lambda^{-1}} \d_\lambda f_\lambda b_2$$
But by assumption $f_{\lambda^{-1}} f_\lambda b_2 = - t_{\lambda^{-1}} f_\lambda b_2 + t_{\lambda^{-1}} \d_\lambda f_\lambda b_2$, so:
\begin{eqnarray*}
- w t_{\lambda^{-1}} f_\lambda b_2 + w t_{\lambda^{-1}} \d_\lambda f_\lambda b_2 & = & - f_\lambda b_2 - \d_\lambda f_\lambda b_2
- w t_{\lambda^{-1}} f_\lambda b_2 + w t_{\lambda^{-1}} \d_\lambda f_\lambda b_2\\
& & - f_\lambda t_{\lambda^{-1}} f_\lambda b_2 + f_\lambda t_{\lambda^{-1}} \d_\lambda f_\lambda b_2
\end{eqnarray*}
One thus has:
\begin{eqnarray*}
f_\lambda b_2 + \d_\lambda f_\lambda b_2 & = & - t_\lambda f_{\lambda^{-1}} f_\lambda b_2 + t_\lambda f_{\lambda^{-1}} \d_\lambda f_\lambda b_2 \\
& = & f_\lambda b_2 - \d_\lambda f_\lambda b_2 + t_\lambda f_{\lambda^{-1}} \d_\lambda f_\lambda b_2
\end{eqnarray*}
But by Step \ref{st:d2=0quadratique:naturelles} which applies here thanks to the second assumption, one has $f_{\lambda^{-1}} \d_\lambda f_\lambda b_2 = - t_{\lambda^{-1}} \d_\lambda f_\lambda b_2$, so one finds $3 \d_\lambda f_\lambda b_2 = 0$. Since $V$ is $3$-torsion-free, we are done.
\end{proofclaim}

\begin{step}\label{st:d2=0quadratique:Z1leqwZ2}
$Z_1 \leq w \cdot Z_2$; in particular if $a_1 \in Z_1$, then $f_\lambda a_1 = - t_\lambda a_1$.
\end{step}
\begin{proofclaim}
Note that the second claim follows immediately from the first and Step \ref{st:d2=0quadratique:naturelles}.
So let $a_1 \in Z_1$. We apply formula (\ref{eq:d2=0quadratique:formule3}) of Step \ref{st:d2=0quadratique:formules} with $a_2 = a_1$ and $\mu = \lambda$:
$$f_{\lambda^{-1}} f_\lambda a_1 + \d_{\lambda^{-1}} f_{\lambda^{-1}} f_\lambda a_1 = - t_{\lambda^{-1}} f_\lambda a_1$$
or equivalently put $u_{\lambda^{-1}} f_{\lambda^{-1}} f_\lambda a_1 = - t_{\lambda^{-1}} f_\lambda a_1$. It follows that $f_{\lambda^{-1}} f_\lambda a_1 = - u_{-\lambda^{-1}} t_{\lambda^{-1}} f_\lambda a_1 = - t_{\lambda^{-1}} u_{-\lambda} f_\lambda a_1$. Since $f_\lambda a_1 \in Z_2$, one finds:
$$f_{\lambda^{-1}} f_\lambda a_1 = - t_{\lambda^{-1}} f_\lambda a_1 + t_{\lambda^{-1}} \d_\lambda f_\lambda a_1$$

This equation is the first assumption of Step \ref{st:d2=0quadratique:court-circuit}. In order to check the second assumption we go back to formula (\ref{eq:d2=0quadratique:formule2}) of Step \ref{st:d2=0quadratique:formules}, which rewrites as follows:
$$- a_1 + w f_\lambda a_1 - t_{\lambda^{-1}} f_\lambda a_1 + t_{\lambda^{-1}} \d_\lambda f_\lambda a_1 = - w t_\lambda a_1$$
or:
$$(w t_\lambda -1) a_1 + (wt_\lambda -1) t_{\lambda^{-1}} f_\lambda a_1 + t_{\lambda^{-1}} \d_\lambda f_\lambda a_1 = 0$$
We apply $(w t_\lambda + 1)$; there remains:
$$- 2 a_1 - 2 t_{\lambda^{-1}} f_\lambda a_1 + (w t_\lambda + 1) t_{\lambda^{-1}} \d_\lambda f_\lambda a_1 = 0$$
This implies in particular that $\d_\lambda f_\lambda a_1 \in w \cdot Z_2$: which is the second assumption needed to apply Step \ref{st:d2=0quadratique:court-circuit} to $b_2 = a_1$ and $\mu = \lambda$.

So one finds $\d_\lambda f_\lambda a_1 = 0$. This means that $\d_\lambda^2 w a_1 = 0$, and this does not depend on $\lambda$. Let us polarize like in Variation \ref{v:repartitioncarres}, that is let us replace $\lambda$ by $\lambda + \mu$; one finds $2 \d_\lambda \d_\mu w a_1 = 0$. Since $V$ is $2$-torsion-free, one has that for all $\lambda, \mu \in \K^\times$, $\d_\lambda \d_\mu w a_1 = 0$, and therefore $w a_1 \in Z_2$.
\end{proofclaim}

We now finish the proof. Let $a_2 \in Z_2$. Formula (\ref{eq:d2=0quadratique:formule3}) of Step \ref{st:d2=0quadratique:formules} is:
$$- \d_{\mu^{-1}} a_2 + f_{\mu^{-1}} f_\mu a_2 + \d_{\mu^{-1}} f_{\mu^{-1}} f_\mu a_2 = - t_{\mu^{-1}} f_\mu a_2 + t_{\mu^{-1}} f_\mu \d_{\mu^{-1}} a_2$$
But since $\d_{\mu^{-1}} a_2 \in Z_1$, one has by Step \ref{st:d2=0quadratique:Z1leqwZ2} that $t_{\mu^{-1}} f_\mu \d_{\mu^{-1}} a_2 = - \d_{\mu^{-1}} a_2$. So one has:
$$f_{\mu^{-1}} f_\mu a_2 + \d_{\mu^{-1}} f_{\mu^{-1}} f_\mu a_2 = - t_{\mu^{-1}} f_\mu a_2$$
or $u_{\mu^{-1}} f_{\mu^{-1}} f_\mu a_2 = - t_{\mu^{-1}} f_\mu a_2$, so that:
\begin{eqnarray*}
f_{\mu^{-1}} f_\mu a_2 & = & - u_{- \mu^{-1}} t_{\mu^{-1}} f_\mu a_2 = - t_{\mu^{-1}} f_\mu a_2 + \d_{\mu^{-1}} t_{\mu^{-1}} f_\mu a_2\\
& = & - t_{\mu^{-1}} f_\mu a_2 + t_{\mu^{-1}} \d_\mu f_\mu a_2
\end{eqnarray*}
which is the first assumption of Step \ref{st:d2=0quadratique:court-circuit}. To check the second assumption, recall that $\d_\mu f_\mu a_2 \in Z_1 \leq w \cdot Z_2$.
It follows from Step \ref{st:d2=0quadratique:court-circuit} applied to $b_2 = a_2$ that $\d_\mu f_\mu a_2 = 0$, that is $\d_\mu^2 w a_2 = 0$.
Here again one polarizes, replacing $\mu$ by $\lambda + \mu$, and one finds $w \cdot Z_2 \leq Z_2$.

So $Z_2$ is $\<U, w\> = G$-invariant; clearly $G$ centralizes $V/Z_2$, so $i$ does too. But $i$ inverts $V$, and since $V$ is $2$-torsion-free, it follows that $V = Z_2$.
\end{proof}

\begin{remark*}\
\begin{itemize}
\item
The assumption that the characteristic is not $3$ appears twice: after Step \ref{st:d2=0quadratique:CVG=0}, in order to bound the $U$-length by $3$, and in Step \ref{st:d2=0quadratique:court-circuit}.
One may wonder what happens in characteristic $3$.
\item
If $\K$ is finite, the classification of $\SL_2(\K)$-modules (Steinberg's tensor product theorem) should imply that only the sums of copies of the natural representation and of trivial modules meet the assumption.
\item
If $\K$ is infinite, I do not know. One should first study the actions of $\SL_2(\F_3(X))$, and I hope that some knowledgeable reader will find the question interesting.
\end{itemize}
\end{remark*}

However and in spite of the Theme, characteristic $3$ is as far as quadratic actions are concerned a special case.

\section{Towards the Algebra}\label{S:algebra}

\subsection{Algebrica}

\begin{variation}\label{v:algebrica}
Let $\K$ be a field of characteristic $\neq 2$ with more than three elements, $G = \SL_2(\K)$, and $V$ be a simple $G$-module of $U$-length $2$. Then the action of $\SL_2(\K)$ induces an action of $\sl_2(\K)$ on $V$ of $\fu$-length $\leq 2$, meaning that $\fu^2 \cdot V = 0$.
\end{variation}
\begin{proof}
\setcounter{equation}{0}
We shall of course argue directly, without using the Theme. Since $V$ is simple, $V$ is $2$-divisible and $2$-torsion-free; moreover $i$ either centralizes or inverts it.
We work in $\End(V)$.

From the relations $u_\lambda w u_{\lambda^{-1}} w u_\lambda w = t_\lambda$, which may be written $u_\lambda w u_{\lambda^{-1}} = t_\lambda w u_{-\lambda} w$, we derive:
$$w + \d_\lambda w + w \d_{\lambda^{-1}} + \d_\lambda w \d_{\lambda^{-1}} = i t_\lambda - t_\lambda w \d_\lambda w$$
which rewrites as:
\begin{equation}\label{e:algebrica1}
i t_\lambda - w = \d_\lambda w + w \d_{\lambda^{-1}} + \d_\lambda w \d_{\lambda^{-1}} + t_\lambda w \d_\lambda w
\end{equation}
We apply $\d_{\lambda^{-1}}$ to the right:
\begin{equation}\label{e:algebrica2}
(i t_\lambda - w)\d_{\lambda^{-1}} = \d_\lambda w \d_{\lambda^{-1}} + t_\lambda w \d_\lambda w \d_{\lambda^{-1}} = (1 + t_\lambda w) \d_\lambda w \d_{\lambda^{-1}}
\end{equation}
There are two cases.
\begin{itemize}
\item
If $i$ centralizes $V$ then $(t_\lambda w)^2 = 1$ and $(1-t_\lambda w)(1 + t_\lambda w) = 0$, hence:
$$0 = (1 - t_\lambda w) (t_\lambda - w)\d_{\lambda^{-1}} = (t_\lambda - w - w + t_\lambda) \d_{\lambda^{-1}}$$
Dividing by $2$, one finds $t_\lambda \d_{\lambda^{-1}} = w \d_{\lambda^{-1}}$.
We apply $\d_\lambda$ to the left in (\ref{e:algebrica1}):
$$\d_\lambda t_\lambda - \d_\lambda w = \d_\lambda w \d_{\lambda^{-1}} + \d_\lambda t_\lambda w \d_\lambda w = t_\lambda \d_{\lambda^{-1}} w \d_\lambda w = 0$$
It follows that $\d_\lambda w = \d_\lambda t_\lambda = t_\lambda \d_{\lambda^{-1}} = w\d_{\lambda^{-1}}$, or $u_\lambda = w u_{\lambda^{-1}} w$.

Hence $t_\lambda = u_\lambda w u_{\lambda^{-1}} w u_\lambda w = u_\lambda^3 w$, and $u_{3\lambda} = t_\lambda w$ has order dividing $2$; in particular $u_{6\lambda} = 1$. The normal closure of unipotent elements is $G$: so if the characteristic is not $3$ one has $G = \{1\}$ in $\End V$. If the characteristic is $3$ then $t_\lambda w = 1$ and $w = t_\lambda$; in particular $w = t_1 = 1$. But since $\K > \F_3$, the normal closure of $w$ is $G$, which therefore centralizes $V$. In this case, $\sl_2(\K)$ acts trivially.
\item
If $i$ inverts $V$, then $(t_\lambda w)^2 = -1$ and $(1 + t_\lambda w)^2 = 2 t_\lambda w$. One deduces from (\ref{e:algebrica2}):
\begin{eqnarray*}
(1 + t_\lambda w) (-t_\lambda - w) \d_{\lambda^{-1}} & = & 2 t_\lambda w \d_\lambda w \d_{\lambda^{-1}}\\
= - (1 + t_\lambda w)(t_\lambda + w)\d_{\lambda^{-1}}\\
= -(t_\lambda + w + w - t_\lambda)\d_{\lambda^{-1}} &= & - 2 w \d_{\lambda^{-1}}
\end{eqnarray*}
Hence $t_\lambda \d_{\lambda^{-1}} + \d_\lambda w \d_{\lambda^{-1}} = 0$. We go back to (\ref{e:algebrica1}), which rewrites as:
$$- t_\lambda - w = \d_\lambda w + w \d_{\lambda^{-1}} + t_\lambda w \d_\lambda w - t_\lambda \d_{\lambda^{-1}}$$
or $(1 + t_\lambda w) \d_\lambda w + (w - t_\lambda) \d_{\lambda^{-1}} + (t_\lambda + w) = 0$.
We apply $(1 - t_\lambda w)$ to the left:
$2 \d_\lambda w + 2 w \d_{\lambda^{-1}} + 2 t_\lambda = 0$.
\end{itemize}
From now on we suppose $i = -1$, so that:
$$\d_\lambda w + w \d_{\lambda^{-1}} = - t_\lambda$$

With this equation we can reconstruct an action of $\sl_2 (\K)$. Let indeed $x_\lambda = \d_\lambda$, $y_\lambda = w \d_\lambda w$, and $h_\lambda = w \d_\lambda - \d_\lambda w$.
We check that we do get a copy of the Lie ring. Since the $U$-length is $2$ it is clear that $\d_{\lambda + \mu} = \d_\lambda + \d_\mu$: which proves the additivity of the maps $\lambda \mapsto x_\lambda$, $\lambda \mapsto y_\lambda$, and $\lambda \mapsto h_\lambda$. It remains to check the bracket identities. Clearly $[x_\lambda, x_\mu] = [y_\lambda, y_\mu] = [h_\lambda, h_\mu] = 0$.

Now since $\d_\lambda w + w \d_{\lambda^{-1}} = - t_\lambda$, one has in particular:
\begin{eqnarray*}
t_\lambda t_\mu & = & (\d_\lambda w + w \d_{\lambda^{-1}})(\d_\mu w + w \d_{\mu^{-1}})\\
& = & \d_\lambda w \d_\mu w + w \d_{\lambda^{-1}} w \d_{\mu^{-1}}\\
= t_{\lambda\mu} & = & - \d_{\lambda \mu} w - w \d_{(\lambda\mu)^{-1}}
\end{eqnarray*}
so that:
$$(\d_\lambda w \d_\mu + \d_{\lambda\mu}) = w (\d_{(\lambda\mu)^{-1}} + \d_{\lambda^{-1}} w \d_{\mu^{-1}}) w$$

Let $q = \d_\lambda w \d_\mu + \d_{\lambda\mu}$: one thus has $\d q = \d w q = 0$. But since $\d w + w \d = -1$, one has $- q = \d w q + w \d q = 0$, whence $q = 0$, that is $\d_\lambda w \d_\mu = - \d_{\lambda\mu}$.
Hence:
\begin{eqnarray*}
[h_\lambda, x_\mu] & = & (w\d_\lambda - \d_\lambda w) \d_\mu - \d_\mu (w\d_\lambda - \d_\lambda w)\\
& = & - \d_\lambda w \d_\mu - \d_\mu w \d_\lambda\\
& = & 2 \d_{\lambda\mu} = 2 x_{\lambda\mu}
\end{eqnarray*}
The similar verification for $[h_\lambda, y_\mu]$ is not any harder. Finally:
$$[x_\lambda, y_\mu] = \d_\lambda w \d_\mu w - w \d_\mu w \d_\lambda
= - \d_{\lambda\mu} w + w \d_{\lambda\mu}
= h_{\lambda\mu}$$
We do retrieve an action of $\sl_2(\K)$. Clearly $\fu^2 \cdot V = 0$.
\end{proof}

\begin{remark*}
One could have with extra arguments avoided the simplicity assumption; these would have involved a few cohomological computations which look alien to the core of the matter.
What the proof given here really shows, is that turning a $G$-module into a $\fg$-module is likely to be harder than turning a $G$-module into a $\K G$-module.
\end{remark*}

\subsection{Logarithmic Variation}

The following should not be compared to Variation \ref{v:d2=0quadratique}.

\begin{variation}\label{v:sl2x2=0quadratique}
Let $\K$ be a field of characteristic $\neq 2$, $\fg = \sl_2(\K)$, $\fb$ be a Borel subring, and $V$ be a $\fb$-module. Suppose that $x^2 \cdot V = 0$. Then $\fu^2 \cdot V = 0$.
\end{variation}
\begin{proof}
Let $\lambda$ and $\mu$ be in $\K$. Then:
$$x_{\frac{\lambda}{2}} x = [h_{\frac{\lambda}{4}}, x] x = - x h_{\frac{\lambda}{4}} x = - x [h_{\frac{\lambda}{4}}, x] = - x x_{\frac{\lambda}{2}}$$
So $x x_\lambda$ annihilates $V$. Now:
$$x_\lambda x_\mu = [h_{\frac{\lambda}{2}}, x] x_\mu = - x h_{\frac{\lambda}{2}} x_\mu = - x [h_{\frac{\lambda}{2}}, x_\mu] = - x x_{\lambda\mu} = 0$$
which means that $\fu^2 \cdot V = 0$.
\end{proof}

\begin{variation}\label{v:sl2quadratiquecoherente}
Let $\K$ be a field of characteristic $\neq 2$, $\fg = \sl_2(\K)$, and $V$ be a $\fg$-module. Suppose that $x^2 \cdot V = 0$. Then for all $\lambda \in \K^\times$, $\ker x_\lambda = \ker x$ and $\im x_\lambda = \im x$.
\end{variation}
\begin{proof}
By Variation \ref{v:sl2x2=0quadratique}, observe that $\fu^2$ annihilates $V$. Then in $\End V$:
$$x_\lambda = [h_{\frac{\lambda}{2\mu}}, x_\mu] = [[x_{\mu}, y_{\frac{\lambda}{2\mu^2}}], x_\mu] = 2 x_{\mu} y_{\frac{\lambda}{2\mu^2}} x_\mu$$
In particular, $\ker x_\mu \leq \ker x_\lambda$ and $\im x_\lambda \leq \im x_\mu$.
\end{proof}

\begin{variation}\label{v:log2Thema}
Let $\K$ be a field of characteristic $\neq 2, 3$, $\fg = \sl_2(\K)$, and $V$ be a simple $\fg$-module. Suppose that $V$ has $x$-length $2$, meaning that $x^2 \cdot V = 0$. Then there exists a $\K$-vector space structure on $V$ making it isomorphic to $\Nat \sl_2(\K)$.
\end{variation}
\begin{proof}
The proof starts here. By simplicity, $\Ann_V (\fg) = 0$; by our assumptions on the base field, $V$ is $6$-torsion-free.
\setcounter{equation}{0}

\begin{step}
$hx = x$ and $(h-1)h(h+1)=0$.
\end{step}
\begin{proofclaim}
One proves by induction in the enveloping ring:
\begin{equation*}
y^i x = xy^i - i(h+i-1) y^{i-1}
\end{equation*}
This equation holds for $i = 0$; one deduces:
\begin{equation}\label{eq:sl2quadratiqueprincipale}
y^i x^2 = x^2y^i - 2i(h+i-2) xy^{i-1} + i(i-1)(h+i-1)(h+i-2)y^{i-2}
\end{equation}
which holds for all $i \geq 0$. We take $i = 1$ in (\ref{eq:sl2quadratiqueprincipale}); one finds $0 = 0 - 2 (h-1) x$, and since $V$ is $2$-torsion-free:
\begin{equation}\label{hx}
hx = x
\end{equation}
We now take $i = 2$ in (\ref{eq:sl2quadratiqueprincipale}); one finds $0 = 0 - 4hxy + 2 (h+1)h$, whence by (\ref{hx}),
$2xy = (h+1)h$.
In particular,
$(h-1)h(h+1) = 2 (h-1) xy = 2 (hx -x)y = 0$.
\end{proofclaim}

Here appears the assumption that the characteristic is not $3$. Recall that for $i \in \Z$ one lets $E_i = \{a \in V: h\cdot a = iv\}$.

\begin{step}
$V = E_{-1} \oplus E_1$ and $\ker x = E_1$.
\end{step}
\begin{proofclaim}
By simplicity, $V$ is $2$-divisible and $2$-torsion-free. Since $(h-1)h(h+1)=0$, one has $V = E_{-1} \oplus E_0 \oplus E_1$; the corresponding projectors are respectively $\frac{1}{2}h(h-1)$, $1 - h^2$, and $\frac{1}{2}h(h+1)$.

If $a_0 \in E_0$, one has $x_\lambda \cdot a_0 \in E_2$; since $V$ is $3$-torsion-free, $E_2 = 0$. So $E_0$ is annihilated by $x_\lambda$ and similarly by $y_\mu$: it follows that $E_0 \leq \Ann_V(\fg) = 0$. Hence $V = E_{-1} \oplus E_1$ (the projectors, namely $\frac{1}{2}(1-h)$ and $\frac{1}{2}(1+h)$, still require $V$ to be $2$-divisible).

We see that $E_1 \leq \ker x$; let us prove the converse. Let $a \in \ker x$; let us write $a = a_{-1} + a_1$ with obvious notations. Then $0 = x\cdot a = x \cdot a_{-1}$, so $a_{-1} \in E_{-1} \cap \ker x$. But since $E_{-1} \leq \ker y$, one finds:
\begin{equation*}
- a_{-1} = h \cdot a_{-1} = xy \cdot a_{-1} - yx \cdot a_{-1} = 0
\end{equation*}
hence $a_{-1} = 0$, that is $a \in E_1$.
\end{proofclaim}

\begin{localnotation}
For $\lambda \in \K$ and $v_i \in E_i$, let:
$$\lambda \cdot v_i = i h_\lambda \cdot v_i \in E_i$$
\end{localnotation}

\begin{step}
This defines an action of $\K$ on $V$; $\sl_2(\K)$ is linear.
\end{step}
\begin{proofclaim}
This is clearly additive in $v_i$ and $\lambda$; it therefore suffices to prove multiplicativity in $\lambda$. Let $\lambda, \mu$ in $\K$.

If $a_1 \in E_1$, one has $\lambda \cdot a_1 = h_\lambda \cdot a_1 = xy_\lambda \cdot a_1 = x_\lambda y \cdot a_1$. Hence:
$$\lambda \cdot (\mu \cdot a_1) = h_\lambda h_\mu \cdot a_1 = x_\lambda y x y_\mu \cdot a_1 = - x_\lambda h y_\mu \cdot a_1 = x_\lambda y_\mu \cdot a_1 = h_{\lambda\mu} \cdot a_1 = (\lambda\mu) \cdot a_1$$
Similarly, for $a_{-1} \in E_{-1}$, $\lambda \cdot a_{-1} = -h_\lambda \cdot a_{-1} = y_\lambda x \cdot a_{-1} = y x_\lambda \cdot a_{-1}$, whence:
$$\lambda \cdot (\mu \cdot a_{-1}) = h_\lambda h_\mu \cdot a_{-1} = y_\lambda x y x_\mu \cdot a_{-1} = y_\lambda x_\mu \cdot a_{-1} = - h_{\lambda\mu} \cdot a_{-1} = (\lambda\mu) \cdot a_{-1}$$
and multiplicativity is proved.

We now show that the action of $\sl_2(\K)$ is linear. The linearity of $h_\lambda$ is obvious; so it suffices to prove that of $x$ and $y$. Let $\lambda \in \K$. The linearity of $x$ on $E_1$ is obvious; now if $a_{-1} \in E_{-1}$, one has:
$$\lambda \cdot (x\cdot a_{-1}) = h_\lambda x \cdot a_{-1} = x y_\lambda x \cdot a_{-1} = - x h_\lambda \cdot a_{-1} = x \cdot (\lambda \cdot a_{-1})$$
The linearity of $y$ on $E_{-1}$ is obvious; if $a_1 \in E_1$, one has:
\[\lambda \cdot (y \cdot a_1) = - h_\lambda y \cdot a_1 = y x_\lambda y \cdot a_1 = y h_\lambda \cdot a_1 = y \cdot (\lambda \cdot a_1)\qedhere\]
\end{proofclaim}

This completes the proof.
\end{proof}

\begin{remark*}
One could also directly prove that a suitable action of $\sl_2(\K)$ induces an action of $\SL_2(\K)$; this would be a converse to Variation \ref{v:algebrica}. One would let $u_\lambda = x_\lambda$ and $w = x - y$. We leave the pleasure of details to the reader; the computations are longer than those of Variation \ref{v:log2Thema}, and the point of going to the group in order to study the Lie ring is disputable.
\end{remark*}

\begin{variation}\label{v:sl2quadratique}
Let $\K$ be a field of characteristic $\neq 2, 3$, $\fg = \sl_2(\K)$, and $V$ be a $\fg$-module of $x$-length at most $2$, meaning that $x^2 \cdot V = 0$. If $\K$  has characteristic $0$, suppose in addition that $V$ is $3$-torsion-free.\inmargin{Fixed} Then $V = \ker h \oplus \ker (h-1)(h+1)$ where $\ker h = \Ann_V(\sl_2(\K))$, and there exists a $\K$-vector space structure on $\ker (h-1)(h+1)$ making it isomorphic to a direct sum of copies of $\Nat \sl_2(\K)$. In particular, $\ker x = \ker x_\lambda$ for all $\lambda \in \K^\times$.
\end{variation}
\begin{proof}
Let $\bar{V} = V/\Ann_V(\fg)$. By perfectness, one has $\Ann_{\bar{V}} (\fg) = 0$. One then reads the proof of Variation \ref{v:log2Thema} again, and sees that simplicity was first used in order to kill $\Ann_V(\fg)$ and $6$-torsion, and then in order to guarantee $2$-divisibility. So one still has $E_0(\bar{V}) = 0$ and $2 \bar{V} \leq E_{-1}(\bar{V}) \oplus E_1(\bar{V})$.
In particular if $a_0 \in E_0(V)$ then $\overline{a_0} = 0$, that is $E_0(V) = \Ann_V(\fg)$.

The proof of Variation \ref{v:log2Thema} constructs for all $\bar{a}_1 \in E_1(\bar{V}) \setminus \{0\}$ a $\K$-linear structure on $\< \fg \cdot \bar{a}_1\>$ such that $\sl_2(\K)$ acts naturally; this also works for $\bar{a}_{-1} \in E_{-1}(\bar{V}) \setminus \{0\}$. In particular, $E_{-1}(\bar{V}) \oplus E_1(\bar{V})$ is a direct sum of vector planes, and so is $2$-divisible.
If $\bar{a} \in \bar{V}$, there is therefore $\bar{b} \in \bar{V}$ such that $2\bar{a} = 4\bar{b}$. Since $\bar{V}$ is $2$-torsion-free, $\bar{a} = 2\bar{b} \in 2 \bar{V}$ and $\bar{V} = 2\bar{V} = E_{-1}(\bar{V}) \oplus E_1(\bar{V})$.

We go back up to $V$\inmargin{Revised} and show that $V = E_{-1} (V) \oplus E_0 (V) \oplus E_1(V)$. Let $a_1 \in \pi^{-1}(E_1(\bar{V}))$. Then $h \cdot \overline{a_1} = \overline{a_1}$ so there is $a_0 \in \Ann_V(\fg) = E_0(V)$ such that $h\cdot a_1 = a_1 + a_0$. Hence $a_1 = (a_1 + a_0) - a_0$ with $h\cdot (a_1 + a_0) = a_1 + a_0$, and $a_1 \in E_0(V) + E_1(V)$. Similarly $\pi^{-1}(E_{-1}(\bar{V})) \leq E_{-1}(V) + E_0(V)$.
Hence $V = \pi^{-1}(\bar{V}) = \pi^{-1} (E_{-1}(\bar{V}) \oplus E_1(\bar{V})) \leq E_{-1}(V) + E_0(V) + E_1(V)$.

The latter sum is direct, for if one has a relation $a_{-1} + a_0 + a_1 = 0$ with obvious notations, then applying $h$ twice one finds $a_{-1} + a_1 = -a_{-1} + a_1 = 0$ whence $2 a_1 = 0$. But $2$ is invertible in $\K$ so $2 h_{\frac{1}{2}} \cdot a_1 = 0 = h \cdot a_1 = a_1$ and $a_{-1} = a_0 = 0$ as well.
Hence $V = E_{-1}(V) \oplus E_0 (V) \oplus E_1(V)$.

We also claim that $E_{-1}(V) \oplus E_1(V)$ is $\fg$-invariant. If $a_1 \in E_1(V)$, then write $x \cdot a_1 = b_{-1} + b_0 + b_1$ with obvious notations and apply $h$. One finds $h x \cdot a_1 = 3 x \cdot a_1 = - b_{-1} + b_1 = 3 b_{-1} + 3 b_1$ whence $2b_1 = 4 b_{-1} = 0$, but applying $h_{\frac{1}{2}}$ this results in $b_1 = b_{-1} = 0$. There remains $x \cdot a_1 = b_0$ with $3 b_0 = 0$. If $\K$ has finite characteristic $p \neq 2, 3$, then $b_0 \in \fg \cdot V$ implies $p b_0 = 0$ and therefore $b_0 = 0$. If $\K$ has characteristic $0$ then by assumption on $V$, $b_0 = 0$. In either case $x \cdot a_1 = 0$ and this means that $E_1(V) \leq \ker x$. Notice that we did not use quadraticity of $x$, so $E_{-1}(V) \leq \ker y$ similarly. Hence $E_{-1}(V) \oplus E_1(V)$ is $\<x, y\>$-invariant.
Moreover, using quadraticity of $x$ and Variation \ref{v:sl2quadratiquecoherente}, $E_1(V) \leq \ker x = \ker x_\lambda$, so $E_{-1}(V) \oplus E_1(V)$ is $\<\{x_\lambda : \lambda a\in \K\}, y\> = \fg$-invariant. 

Finally $E_{-1}(V) \oplus E_1(V)$ is a $\fg$-submodule disjoint from $E_0(V) = \Ann_V(\fg)$, so it is isomorphic to $\bar{V}$: it is a direct sum of copies of the natural representation.
\end{proof}

\subsection{Characteristic 3}

\begin{remark*}
As opposed to the Theme to which it is a Lie ring analog, Variation \ref{v:log2Thema} does not hold in characteristic $3$.

Let indeed $\K$ be a field of characteristic $3$.
Let $V = \K e_2 \oplus \K e_0 \oplus \K e_1$; let $x$ and $y$ act by:
$$\left\{
\begin{array}{rcl}
x \cdot e_2 & = & e_1\\
x \cdot e_0 & = & 0\\
x \cdot e_1 & = & 0\\
\end{array}
\right.,\qquad
\left\{
\begin{array}{rcl}
y \cdot e_2 & = & e_0\\
y \cdot e_0 & = & e_1\\
y \cdot e_1 & = & e_2\\
\end{array}
\right.$$
and extend linearly. One may check that this does define an action of $\sl_2(\K)$ where $x^2$ is trivial.
\end{remark*}

\tikzset{auto}
\[\begin{tikzpicture}
 \node (E0) {$E_0$};
 \node (E-1) [left of=E0, below of=E0] {$E_{-1}$};
 \node (E1) [right of=P, below of=E0] {$E_1$};
 \draw[->, bend left] (E-1) to node {$y$} (E0);
 \draw[->, bend left] (E0) to node {$y$} (E1);
 \draw[->] (E-1.8) to node {$x$} (E1.170);
 \draw[<-, bend right] (E-1.350) to node [swap] {$y$} (E1.195);
\end{tikzpicture}
\qedhere\]

One will in particular note that $x^2 = 0 \neq y^2$: \emph{this representation of the Lie ring cannot come from a representation of the group}.

\begin{variation}\label{v:sl2car3structurepartielle}
Let $\K$ be a field of characteristic $3$, $\fg = \sl_2(\K)$, and $V$ be a simple $\fg$-module with $x^2 = 0$ in $\End V$.
Then $E_{-1} \oplus E_1$ may be equipped with a $\K$-vector space structure such that, saying that $\K$ annihilates $E_0$, the maps $h_\lambda$ and $x_\lambda$ are everywhere linear (the $y_\lambda$'s a priori only on $E_1$).
\end{variation}
\begin{proof}
We go back to the proof of Variation \ref{v:log2Thema}; in characteristic $3$ one still has the equations $(h-1)h(h+1) = 0$ and $hx = x$. $V$ being $2$-divisible (it has exponent $3$), it follows that $V = E_{-1} \oplus E_0 \oplus E_1$, and $x\cdot V \leq E_1$. In particular, $x \cdot E_0 \leq E_{-1}\cap E_1 = 0$, and $x\cdot E_1 \leq E_0 \cap E_1 = 0$. This proves $E_0 \oplus E_1 \leq \ker x$.

Now suppose $a_{-1} \in E_{-1} \cap \ker x$. Then:
$$-a_{-1} = h\cdot a_{-1} = xy \cdot a_{-1} - yx \cdot a_{-1}$$
Since $y \cdot a_{-1} \in E_1 \leq \ker x$, one finds $a_{-1} = 0$: hence $\ker x = E_0 \oplus E_1$.

Therefore the module is as in the diagram above. On $E_{-1} \oplus E_1$ one defines the same linear structure as in Variation \ref{v:log2Thema}: this still makes sense as one will check.
\end{proof}

\begin{remark*}
One can't go any further. Let indeed $\K > \F_3$ be a field of characteristic $3$ and take three copies of $\K^3$, denoted $E_i$, the elements of which are the $\lambda_i$'s for $\lambda \in \K$, $i \in \{-1, 0, 1\}$; one identifies $0_{-1} = 0_0 = 0_1$.

Let $\sigma$ be an additive map from $\K$ to $\K$. We then define an action of $\sl_2(\K)$ as follows:
$$\left\{\begin{matrix}
x_\lambda \cdot (\mu_1) & = & 0\\
x_\lambda \cdot (\mu_0) & = & 0\\
x_\lambda \cdot (\mu_{-1}) & = & (\lambda\mu)_1
\end{matrix}\right.
,\quad
\left\{\begin{matrix}
y_\lambda \cdot (\mu_1) & = & (\lambda\mu)_{-1}\\
y_\lambda \cdot (\mu_0) & = & (\lambda\mu)_1\\
y_\lambda \cdot (\mu_{-1}) & = & (\sigma(\lambda\mu))_0
\end{matrix}\right.$$
Since $\sigma$ is additive, this does define a $\fg$-module where $x^2 = 0$. One can actually make $V$ simple by taking $\sigma$ to be injective; in general, starting with any element of $\im \sigma$, one can reconstruct $E_{-1} \oplus (\im \sigma)_0 \oplus E_1$.

If there were a compatible linear structure, $y^3$ would be linear; yet $(y^3)_{|E_0} = \sigma$.
One can chose $\sigma$ so that $\ker (\sigma - \Id)$ has exactly $3$ elements: $\sigma$ will then be linear for no $\K$-vector space structure.

We have just constructed a representation of the Lie ring $\sl_2(\K)$ which cannot come from a representation of the Lie algebra.

There is slightly worse. We now take $\sigma$ to be an additive map such that the cardinal of $\im \sigma$ is strictly less than that of $\K$ (this is possible be $\K$ finite or infinite). One then obtains a simple $\sl_2(\K)$-module of the form $E_{-1} \oplus (\im \sigma)_0 \oplus E_1$. For cardinality reasons, the null weight subgroup cannot be equipped with any $\K$-vector space structure: this explains our embarrassment on $E_0$ in Variation \ref{v:sl2car3structurepartielle}.
\end{remark*}

\begin{remark*}
Observe however that even in characteristic $3$, if \emph{both} $x^2$ and $y^2$ are zero on the simple $\sl_2(\K)$-module $V$, then both $x$ and $y$ annihilate $E_0$. As a consequence and by Variation \ref{v:sl2quadratiquecoherente}, $E_0 \leq \Ann_V(\sl_2(\K)) = 0$. So there exists a $\K$-vector space structure on $V$ making it isomorphic to $\Nat \sl_2(\K)$.
\end{remark*}

One may remove simplicity.\inmargin{From there on, added.}

\begin{variation}\label{v:sl2car3biquadratique}
Let $\K$ be a field of characteristic $3$, $\fg = \sl_2(\K)$, and $V$ be a $\fg$-module with $x^2 = y^2 = 0$ in $\End V$. Then $V = \Ann_V(\fg) \oplus \fg \cdot V$, and there exists a $\K$-vector space structure on $\fg\cdot V$ making it isomorphic to a direct sum of copies of $\Nat \sl_2(\K)$.
\end{variation}
\begin{proof}
We shall first work with $\F_3$, the field with three elements. Let $\fg_1 = \sl_2(\F_3)$ as a Lie subring of $\fg$ and consider the $\fg_1$-module $V$. The $\fg$-analysis will be made in the end.

$V$ need not have exponent $3$. If one reads the computations of Variation \ref{v:log2Thema} again, one will merely expect $2hx = 2x$ and $2(h-1)h(h+1) = 0$.
However $\fg_1 \cdot V$ does have exponent $3$, and so does the ideal generated by $\fg_1$ in $\End V$. In particular one has $(h-1)h(h+1) = 0$ and $hx = x$ in $\End V$; by the quadraticity assumption on $y$ one has $hy = -y$ as well.

Let $\bar{V} = V/\Ann_V \fg_1$. By perfectness of $\fg_1$, $\Ann_{\bar{V}} \fg_1 = 0$. So $3 \bar{V} \leq \Ann_{\bar{V}} \fg_1 = 0$ and $\bar{V}$ has exponent $3$. Of course in $\End \bar{V}$ the equations $(h-1)h(h+1) = 0$, $hx = x$, and $hy = -y$ still hold.

Since $\bar{V}$ is a vector space over $\F_3$ one derives $\bar{V} = E_{-1}(\bar{V}) \oplus E_0(\bar{V}) \oplus E_1(\bar{V})$. But then $x \cdot (E_0(\bar{V}) \oplus E_1(\bar{V})) \leq (E_{-1}(\bar{V}) \oplus E_0(\bar{V})) \cap E_1(\bar{V}) = 0$. Symmetrically, $y$ annihilates $E_{-1}(\bar{V}) \oplus E_0(\bar{V})$. It follows that $E_0(\bar{V}) \leq \ker x \cap \ker y = \Ann_{\bar{V}}\fg_1 = 0$. Therefore $\bar{V} = E_{-1}(\bar{V}) \oplus E_1(\bar{V})$.

As said, $x$ annihilates $E_1(\bar{V})$ and $y$ annihilates $E_{-1}(\bar{V})$. Moreover $x$ is injective on $E_{-1}(\bar{V})$ since for $\overline{a_{-1}} \in E_{-1}(\bar{V}) \cap \ker x$ one has $-\overline{a_{-1}} = h \cdot \overline{a_{-1}} = (xy-yx) \cdot \overline{a_{-1}} = 0$. At this point it is clear that $\bar{V} = E_{-1}(\bar{V}) \oplus E_1(\bar{V})$ is a direct sum of copies of $\Nat \fg_1$.

We go back up to $V$ exactly like in Variation \ref{v:sl2quadratique} and show that $V = E_{-1}(V) \oplus E_0(V) \oplus E_1(V)$.
We also claim that $E_{-1}(V) \oplus E_1(V)$ is $\fg_1$-invariant.
If $a_1 \in E_1(V)$ then a priori using the same notations as in Variation \ref{v:sl2quadratique} one should find $x \cdot a_1 = b_0$ with $b_0 \in \Ann_V\fg_1 = E_0(V)$ of order $3$. But quadraticity of $x$ proved that in $\End V$, $2hx = 2x$. Hence $0 = 2h \cdot b_0 = 2h x \cdot a_1 = 2x \cdot a_1 = 2b_0$. There remains $b_0 = 3 b_0 - 2 b_0 = 0$, and $E_1(V) \leq \ker x$. But since we have assumed that $y$ is quadratic as well, one also has $E_{-1}(V) \leq \ker y$, and this proves that $E_{-1}(V) \oplus E_1(V)$ is $\fg_1$-invariant.

It is now clear that $\fg_1 \cdot V = E_{-1}(V) \oplus E_1(V) \simeq V/E_0(V) \simeq \bar{V}$ as a $\fg_1$-module is a direct sum of copies of $\Nat \fg_1$, and $V = \Ann_V \fg_1 \oplus \fg_1 \cdot V$.

We move to another set of ideas. By Variation \ref{v:sl2quadratiquecoherente}, $\im x = \im x_\lambda$ and $\ker x = \ker x_\lambda$ for all $\lambda \in \K^\times$, and similarly with $y$ and $y_\lambda$. So as a matter of fact, $\Ann_V \fg_1 = \ker x \cap \ker y = \Ann_V \fg$ and $\fg_1 \cdot V = \im x + \im y = \fg \cdot V$.

The same linear construction as in Variation \ref{v:log2Thema} will then provide a suitable $\K$-vector space structure on $\fg \cdot V = E_{-1}(V) \oplus E_1(V)$.
\end{proof}

\medskip
\hrule
\medskip
\noindent
Future variations will explore the symmetric powers of $\Nat \fG$.

\bibliographystyle{plain}
\bibliography{Variationen}

\end{document}